\documentclass{amsart}

\theoremstyle{plain}
\newtheorem{Theorem}{Theorem}
\newtheorem{Corollary}[Theorem]{Corollary}
\newtheorem{Lemma}[Theorem]{Lemma}
\newtheorem{Proposition}[Theorem]{Proposition}

\theoremstyle{definition}
\newtheorem*{definition}{Definition}

\theoremstyle{remark}
\newtheorem*{Remark}{Remark}

\DeclareMathOperator{\Bs}{Bs}

\DeclareMathOperator{\Supp}{Supp}
\DeclareMathOperator{\mult}{mult}

\newcommand{\Z}{\ensuremath{\mathbb{Z}\,}}
\newcommand{\Q}{\ensuremath{\mathbb{Q}\,}}

\newcommand{\C}{\ensuremath{\mathbb{C}\,}}

\newcommand{\Ocal}{\mathcal{O}}

\newcommand{\eqdef}{\overset{\text{def}}{=}}
\newcommand{\fne}{\underset{f}{\equiv}}

\newcommand{\ulc}{\lceil}
\newcommand{\urc}{\rceil}
\newcommand{\llc}{\lfloor}
\newcommand{\lrc}{\rfloor}

\newcommand{\rup}[1]{\ulc #1 \urc}	
\newcommand{\rdn}[1]{\llc #1 \lrc}	

\newcommand{\bskm}{\Bs |K_Y+ \ulc M \urc|}

\begin{document}



\title[Normal surface singularities] 
	{Kawachi's invariant for normal surface singularities}

\author{Vladimir Ma\c{s}ek}
\address{Department of Mathematics, Box 1146, Washington University,
		St. Louis, MO 63130}
\email{vmasek@math.wustl.edu}


\begin{abstract}
  We study a useful numerical invariant of normal surface singularities,
  introduced recently by T. Kawachi. Using this invariant, we give a quick 
  proof of the (well-known) fact that all log-canonical surface singularities
  are either elliptic Gorenstein or rational (without assuming a priori that
  they are \Q-Gorenstein).
  
  In \S 2 we prove effective results (stated in terms of Kawachi's invariant)
  regarding global generation of adjoint linear systems on normal surfaces with
  boundary. Such results can be used in proving effective estimates for global
  generation on singular threefolds.  The theorem of Ein--Lazarsfeld and
  Kawamata, which says that the minimal center of log-canonical singularities
  is always normal, explains why the results proved here are relevant in that
  situation.
\end{abstract}

\maketitle




\setcounter{section}{-1}
\subsection*{Contents}

\begin{enumerate}
  \item[0.] Introduction 
  \item[1.] Kawachi's invariant and log-canonical singularities
  \item[2.] A theorem of Reider type on normal surfaces with boundary
\end{enumerate}

\subsection*{Notations}

\begin{tabbing}
99\=9999999999\=9999999999999999999999999999\kill
  \>$\ulc \cdot \urc$ \> round-up               \\
  \>$\llc \cdot \lrc$ \> round-down             \\
  \>$\{ \cdot \}$     \> fractional part        \\
  \>$f^{-1}D$         \> strict transform (proper transform) \\
  \>$f^*D$            \> pull-back (total inverse image) \\
  \>$\,\equiv$        \> numerical equivalence  \\
  \>$\,\sim$          \> linear equivalence     \\
\end{tabbing}


\section{Introduction}

Let $Y$ be a normal algebraic surface over an algebraically closed field
of arbitrary characteristic. Let $y \in Y$ be a fixed point on $Y$. Let
$f : X \to (Y,y)$ be the minimal resolution of the germ $(Y,y)$ if $y$ is
singular, resp. the blowing-up of $Y$ at $y$ if $y$ is smooth. Kawachi
(\cite{kawachi1}) introduced the following numerical invariant of $(Y,y)$:

\begin{definition}
  $\delta_y = -(Z-\Delta)^2$, where $Z$ is the fundamental cycle of $y$ and 
  $\Delta = f^*K_Y - K_X$ is the canonical cycle (or antidiscrepancy) of $y$.
\end{definition}

In \S 1 we recall several definitions (including the fundamental and the
canonical cycle, etc.); then we study Kawachi's invariant and we give a very
short proof of the fact (well-known to the experts) that ``numerically''
log-canonical surface singularities are automatically \Q-Gorenstein.

In \S 2 we prove several criteria for global generation of linear systems of
the form $|K_Y+\rup{M}|$, $M$ a \Q-divisor on $Y$ such that $K_Y+\rup{M}$ is
Cartier. This type of results was the original motivation for introducing the
invariant $\delta_y$ (see \cite{kawachi1, kawachi2, kawachim}). The main
interest in results of Reider type for \Q-divisors on normal surfaces comes
from work related to Fujita's conjecture for (log-) terminal threefolds, cf.
\cite{elm, matsushita}.  Using the criterion proved in \S 2, together with
other recent results, we can significantly improve the main theorem of
\cite{elm}; we will do so in a forthcoming paper.

The author has benefitted from numerous discussions with L.~B\u{a}descu,
L.~Ein, T.~Kawachi, R.~Lazarsfeld, and N.~Mohan Kumar.


\section{Kawachi's invariant and log-canonical singularities}

In this section we recall several standard definitions and facts regarding
normal surface singularities and we prove a number of elementary lemmas
involving Kawachi's invariant. For convenience, we use M. Reid's recent notes 
\cite{reid} as our main reference.

\vspace{6pt}

{\bf (1.1)}
Let $Y$ be a complete normal algebraic surface (in any characteristic), and
let $f : X \to Y$ be a resolution of singularities of $Y$. We use Mumford's
\Q-valued pullback and intersection theory on $Y$ (cf. \cite{mumford}): if $D$
is any Weil (or \Q-Weil) divisor on $Y$, then $f^*D = f^{-1}D+D_{\text{exc}}$,
where $f^{-1}D$ is the strict transform of $D$ and $D_{\text{exc}}$ is the
(unique) $f$-exceptional \Q-divisor on $X$ such that $f^*D \cdot F = 0$ for
every $f$-exceptional curve $F \subset X$. (The existence and uniqueness of
$D_{\text{exc}}$ follow from the negative definiteness of the intersection
form on exceptional curves.) If $D_1, D_2$ are two \Q-Weil divisors on $Y$,
then $D_1 \cdot D_2 \eqdef f^*D_1 \cdot f^*D_2$. See, e.g., \cite[\S 1]{elm}
for a quick review of this theory. [Here are a few key points: $D \geq 0
\implies f^*D \geq 0$; $D_1 \cdot D_2$ is independent of resolution; if $C$ is
a \Q-divisor on $X$, then $C \cdot f^*D = f_*C \cdot D$; if $D$ is \Q-Cartier,
then the definition of $f^*D$ coincides with the usual one.]

\begin{definition}
   (a) Let $M$ be a \Q-Weil divisor on $Y$. Then $M$ is \emph{nef} if 
         $M \cdot C \geq 0$ for all irreducible curves $C \subset Y$.
         (Equivalently, $M$ is nef if and only if $f^*M$ is nef on the smooth
         surface $X$.)

   (b) Assume that $M$ is nef. Then $M$ is \emph{big} if in addition $M^2>0$
   	 (i.e., if $f^*M$ is big on $X$).
\end{definition}

\vspace{6pt}

{\bf (1.2)}
Now let $y \in Y$ be a fixed point, and let $f : X \to (Y,y)$ be the
\emph{minimal} resolution of the germ $(Y,y)$ if $y$ is singular, resp. the
blowing-up of $Y$ at $y$ if $y$ is smooth. Let $f^{-1}(y) = \cup_{j=1}^{N} F_j$
(set-theoretically); $N = 1$ if $y$ is smooth.

\begin{definition}
   The \emph{fundamental cycle} of $(Y,y)$ is the smallest nonzero effective
   divisor $Z = \sum z_j F_j$ on $X$ (with $z_j \in \Z$) such that 
   $Z \cdot F_j\leq 0,\forall j$ (cf. \cite[4.5]{reid} or \cite[p.132]{artin}).
\end{definition}

Note that $z_j \geq 1, \forall j$, because $\cup F_j$ is connected.

If $y$ is smooth, then $Z = F_1$ ($F_1$ is a $(-1)$-curve in this case).

Let $p_a(Z) = \frac{1}{2} Z \cdot (Z+K_X) + 1$; then $p_a(Z) \geq 0$, cf.
\cite[Theorem 3]{artin}.

\vspace{6pt}

{\bf (1.3)}
Let $K_X$ be a canonical divisor on $X$; then $f_*K_X$ is a canonical divisor
$K_Y$ on $Y$, and $\Delta \eqdef f^*K_Y - K_X$ is an $f$-exceptional \Q-divisor
on $X$, $\Delta = \sum a_jF_j$.
Note that  $\Delta$ is $f$-numerically equivalent to $K_X$
($\Delta \fne -K_X$), i.e. $\Delta \cdot F_j = -K_X \cdot F_j, \forall j$;
in particular, $p_a(Z) = \frac{1}{2} Z \cdot (Z-\Delta) + 1$, or
$Z \cdot (Z-\Delta) = 2 p_a(Z) - 2$.

If $y$ is smooth, then $\Delta = -F_1$. On the other hand, if $y$ is singular
and $f$ is the minimal resolution of $(Y,y)$, then
$\Delta \cdot F_j = -K_X \cdot F_j \leq 0, \forall j$, so that $\Delta$ is
effective (see, e.g., \cite[Lemma 1.4]{elm}). In fact, $\Delta = 0$ if and only
if $y$ is a canonical singularity (= $RDP$, = du Val singularity), cf. the last
part of (1.4) below, and if $y$ is not canonical, then $a_j > 0$ for \emph{all}
$j = 1,\dots ,N$ (again, because $\cup F_j$ is connected).

\begin{definition}
 $\Delta$ is the \emph{canonical cycle} (or \emph{antidiscrepancy}) of $(Y,y)$.
\end{definition}

$\Delta$ is uniquely defined, even though $K_X$ (and, accordingly, $K_Y$)
is defined only up to linear equivalence. This follows from
$\Delta \cdot F_j = -K_X \cdot F_j, \forall j$, and the negative definiteness
of $\Vert F_i \cdot F_j \Vert$.

Assume that $y$ is Gorenstein and non-canonical. Then $\Delta$ has integer
coefficients and $\Delta > 0$; by the definition of the fundamental cycle, 
$\Delta \geq Z$. On the other hand, if $(Y,y)$ is a log-canonical singularity, 
then $\Delta \leq Z$ (see below).

\vspace{6pt}

{\bf (1.4)}
{\bf Definition.}  $y$ is a \emph{rational} singularity if $R^1f_*\Ocal_X = 0$,
or equivalently (cf. \cite[Theorem 3]{artin}), if $p_a(Z)=0$.

\vspace{4pt}

If $y$ is any normal singularity and $D$ is a Cartier divisor on $Y$, then
$f^*D$ has integer coefficients. If $y$ is a rational singularity, then the
converse is also true: if $f^*D$ has integer coefficients, then $D$ is Cartier.
(\emph{Proof:} $f^*D \cdot F_j = 0, \forall j$, and therefore $f^*D$ is trivial
in an open neighborhood of $f^{-1}(y)$, cf. \cite[Lemma 5]{artin}; thus $D$ is
trivial in a punctured open neighborhood of $y$ in $Y$, and therefore $D$ is
Cartier, because $y$ is normal.) In particular, a rational surface singularity
is always \Q-factorial.

If $y$ is a rational singularity, then $\mult_y Y =-Z^2$
(cf. \cite[Corollary 6]{artin}, or \cite[4.17]{reid}); in particular, 
$\text{$y$ is a $RDP$} \iff Z^2 = -2 \iff Z \cdot \Delta = 0 \iff
\Delta \cdot F_j = 0, \forall j \iff \Delta = 0$. (Therefore $RDP \implies$
Gorenstein, because $f^*K_Y = K_X$ has integer coefficients, and therefore
$K_Y$ is Cartier. The converse is also true: rational Gorenstein
$\implies RDP$. Indeed, if $y$ is Gorenstein, then either $\Delta = 0$ or
$\Delta \geq Z$; but $\Delta \geq Z \implies Z \cdot (Z-\Delta) \geq 0
\implies p_a(Z) \geq 1$.)

\vspace{6pt}

{\bf (1.5)} 
{\bf Definition.}  $y$ is an \emph{elliptic} singularity if $R^1f_* \Ocal_X$
is 1-dimensional. It is \emph{elliptic Gorenstein} if in addition it is
Gorenstein.

\begin{Lemma}
  $y$ is elliptic Gorenstein if and only if $Z = \Delta$. 
\end{Lemma}

\begin{proof}
  $\implies$ is proved in \cite[4.21]{reid}. Conversely, if $Z = \Delta$, then
  $K_X+Z = f^*K_Y$, so that $(K_X+Z) \cdot F_j = 0, \forall j$; the proof
  proceeds as in \cite[4.21]{reid} (go directly to Step 3 there).
\end{proof}

{\bf (1.6)}
{\bf Definition.} $y$ is \emph{log-terminal} (resp. \emph{log-canonical}) if
$a_j < 1$ (resp. $a_j \leq 1$), $\forall j$ -- where $\Delta = f^*K_Y-K_X =
\sum a_jF_j$, as before.

In dimension three or higher, one must assume that $y$ is \Q-Gorenstein before
defining $\Delta$ (and log-terminal singularities, etc.) Using Mumford's
definition of $f^*K_Y$, we don't need to make this assumption in the
two-dimensional case. What's more, we will see in a moment that log-canonical
(in our sense) automatically implies \Q-Gorenstein. (This is also clear from
the complete list of all log-canonical singularities: the arguments in 
\cite{alexeev} do not use the \Q-Gorenstein condition.)

\vspace{6pt}

{\bf (1.7)}
{\bf Definition.}  $\delta_y = -(Z-\Delta)^2$.

\vspace{4pt}

Thus $\delta_y \in \Q, \delta_y \geq 0$, and $\delta_y = 0 \iff Z = \Delta
\iff y$ is elliptic Gorenstein, by Lemma~1.

\begin{Lemma}  { \ \ \ }
  \begin{enumerate}
    \item[(a)] $\delta_y = 2-2p_a(Z)-\sum_{j=1}^N (z_j-a_j) K_X \cdot F_j$.
    \item[(b)] $\delta_y = 2-2p_a(Z)+\sum_{j=1}^N a_j (Z-\Delta) \cdot F_j$.
  \end{enumerate}
\end{Lemma}

\begin{proof}
  (a) $\delta_y = -(Z-\Delta)^2 = -Z \cdot (Z-\Delta) + \Delta \cdot (Z-\Delta)
  = 2 - 2p_a(Z) - K_X \cdot (Z-\Delta)
  = 2-2p_a(Z)-\sum_{j=1}^{N}(z_j-a_j) K_X\cdot F_j$ (because $\Delta\fne K_X$).

  (b) is similar.
\end{proof}

\begin{Lemma} \textup{(cf. \cite[Theorem 1]{kawachim}) }
  \begin{enumerate}
    \item[(a)] $\delta_y = 4$ if $y$ is smooth;
    \item[(b)] $\delta_y = 2$ if $y$ is a $RDP$;
    \item[(c)] $0<\delta_y<2$ if $y$ is log-terminal but not smooth or a $RDP$;
    \item[(d)] $0 \leq \delta_y \leq 2$ if $y$ is log-canonical but not smooth.
 \end{enumerate}
\end{Lemma}

\begin{proof}
  (a) $y$ smooth $\implies Z-\Delta=2F_1$, and $F_1^2=-1$; thus $\delta_y=4$. 

  (b) $y$ a $RDP \implies \Delta=0$, and $\delta_y = -Z^2 = \mult_y Y = 2$.

  (c), (d) $\delta_y$ is always $\geq 0$. If $\delta_y = 0$, then $\Delta = Z$,
  so that $a_j = z_j \geq 1, \forall j$; therefore $y$ log-terminal $\implies
  \delta_y > 0$.

  On the other hand, log-canonical $\implies a_j \leq 1 \leq z_j, \forall j$,
  so that $\delta_y \leq 2$ by Lemma~2 above ($K_X\cdot F_j\geq 0,\forall j$,
  because $f$ is the \emph{minimal} resolution). 

  Finally, if $y$ is log-terminal, then $a_j < z_j, \forall j$; then
  $\delta_y = 2 \implies p_a(Z) = 0$ and $K_X \cdot F_j = 0, \forall j$
  (recall that $p_a(Z) \geq 0$ for all normal singularities). Thus $y$ is a
  rational singularity, and the last paragraph of (1.4) shows that $y$ is a 
  $RDP$.
\end{proof}

{\bf Corollary of the proof.} \emph{
  If $y$ is log-terminal then it is rational. If $y$ is log-canonical then it
  is either elliptic Gorenstein or rational. In particular, log-canonical
  implies \Q-Gorenstein.}  (We noted already in (1.4) that rational implies
  \Q-factorial.)

Indeed, the proof above shows that $y$ log-canonical $\implies p_a(Z) \leq 1$
(because $\delta_y \geq 0, z_j-a_j \geq 0, \text{ and } K_X \cdot F_j \geq 0,
\forall j$); moreover, $p_a(Z) = 1 \implies \delta_y = 0 \implies Z = \Delta
\implies y$ is elliptic Gorenstein -- and this can happen only in the
log-canonical, not in the log-terminal case. In all other cases, $p_a(Z) = 0$,
and therefore $y$ is rational.

\vspace{6pt}

{\bf (1.8)}
The invariant $-\Delta^2$ has been considered before; for example, Sakai
\cite{sakai} proved results of Reider type on normal surfaces, using this
invariant. In a sense, $-\Delta^2$ may be viewed as a local analogue of the
Chern number $c_1^2$.

\vspace{4pt}

At least in the rational case, $-\Delta^2$ is closely related to $\mult_y Y$,
via Kawachi's invariant:

\begin{Lemma}
  $-\Delta^2 = -Z^2 + \delta_y + 4(p_a(Z)-1)$
\end{Lemma}
\begin{proof}
  $\Delta = Z - (Z-\Delta)$; therefore  \newline
  $-\Delta^2 = -Z^2 -(Z-\Delta)^2 +2 Z \cdot (Z-\Delta) = 
  -Z^2 + \delta_y + 4(p_a(Z)-1)$
\end{proof}

\begin{Corollary}
  If $y$ is rational, then $p_a(Z) = 0$ and $-Z^2 = \mult_y Y$, so that 
    $$ -\Delta^2 = \mult_y Y - (4 - \delta_y).  $$
\end{Corollary}
In particular, if $y$ is log-terminal but not canonical, we get:
    $$ \mult_y Y - 4 < -\Delta^2 < \mult_y Y - 2.  $$

\vspace{6pt}

{\bf (1.9)}
Since log-canonical surface singularities are classified (see, for example,
\cite{alexeev, crepant}), one could compute $\delta_y$ explicitly in all cases.
Indeed, this was Kawachi's original proof of part (c) in Lemma 3; see
\cite{kawachi1} for the complete list in the log-terminal case.

The computation of $Z, \Delta, \text{ and } \delta_y$ is an easy exercise
in linear algebra. For illustration, we give here the values of $Z, \Delta,
\text{ and } \delta_y$ for the ``truly'' log-canonical (i.e. non-log-terminal)
singularities. The list of all such singularities can be found in 
\cite[p.58]{alexeev}. To simplify notation, we assume that the chains of smooth
rational curves shown in \cite{alexeev} consist of just one curve each. In
each case, we show the dual graph of $\cup F_j$, indicating the
self-intersection numbers $F_j^2$.

Notice that in the ``truly'' log-canonical case $\delta_y$ takes only the
values 0, 1, and 2 (always an integer), and the value of $\delta_y$
distinguishes three types of log-canonical singularities
(with one interesting exception, noted below):

\begin{picture}(290,50)(0,0)   
    \put(0,25){\makebox(0,0)[l]{Notation:}}
    \put(51,25){\circle{8}}
    \put(51,25){\circle{3}}
    \put(59,25){\makebox(0,0)[l]{= smooth elliptic curve; }}
    \put(174,25){\circle{8}}
    \put(182,25){\makebox(0,0)[l]{= smooth rational curve.}}
\end{picture}

\emph{Type 1:} elliptic Gorenstein (cases (4) and (5) in \cite{alexeev})

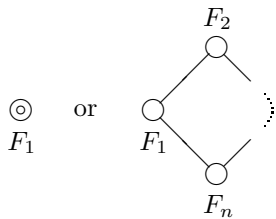
\begin{center}
\begin{picture}(200,85)(0,10)   
    \put(50,50){\circle{8}}
    \put(50,50){\circle{3}}
    \put(45,35){$F_1$}
    \put(70,50){\makebox(0,0)[l]{or}}
    \put(100,50){\circle{8}}
    \put(95,35){$F_1$}
    \put(124,26){\circle{8}}
    \put(119,11){$F_n$}
    \put(124,74){\circle{8}}
    \put(119,82){$F_2$}
    \put(102.82,52.82){\line(1,1){18.56}}
    \put(102.82,47.18){\line(1,-1){18.56}}
    \put(126.82,71.18){\line(1,-1){10}}
    \put(126.82,28.82){\line(1,1){10}}
    \qbezier[7](142,56)(148,50)(142,44)
\end{picture}
\end{center}

$$  Z = \Delta = F_1, \text{ resp. } F_1+\cdots+F_n; \qquad \delta_y = 0.  $$

\vspace{15pt}

\emph{Type 2:}  (case (6) in \cite{alexeev})

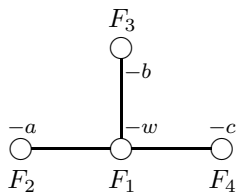
\begin{center}
\begin{picture}(200,80)(0,18)   
    \put(50,40){\circle{8}}
        \put(45,25){$F_2$}
        \put(45,46){\Small $-a$}
    \put(88,40){\circle{8}}
        \put(83,25){$F_1$}
        \put(89,46){\Small $-w$}
    \put(126,40){\circle{8}}
        \put(121,25){$F_4$}
        \put(121,46){\Small $-c$}
    \put(88,78){\circle{8}}
        \put(83,86){$F_3$}
        \put(89,67){\Small $-b$}
    \put(54,40){\line(1,0){30}}
    \put(92,40){\line(1,0){30}}
    \put(88,44){\line(0,1){30}}
\end{picture}
\end{center}
\begin{gather*}
	w \geq 2;\quad (a,b,c)=(3,3,3),\,(2,2,4),\,\text{ or } (2,3,6)	\\
	Z = F_1+(F_2+F_3+F_4) \text{ if } w \geq 3,
	Z = 2F_1+(F_2+F_3+F_4) \text{ if } w = 2;			\\
	\Delta = F_1 + (1-\tfrac{1}{a})F_2 + 
	(1-\tfrac{1}{b})F_3 + (1-\tfrac{1}{c})F_4; \qquad
	\delta_y = 1.
\end{gather*}

\vspace{15pt}

\emph{Type 3:} (cases (7) -- (8) in \cite{alexeev})

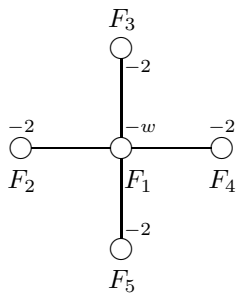
\begin{center}
\begin{picture}(200,120)(0,-20)   
    \put(50,40){\circle{8}}
        \put(45,25){$F_2$}
        \put(45,46){\scriptsize $-2$}
    \put(88,40){\circle{8}}
        \put(89,25){$F_1$}
        \put(89,46){\scriptsize $-w$}
    \put(126,40){\circle{8}}
        \put(121,25){$F_4$}
        \put(121,46){\scriptsize $-2$}
    \put(88,78){\circle{8}}
        \put(83,86){$F_3$}
        \put(89,69){\scriptsize $-2$}
    \put(88,2){\circle{8}}
        \put(83,-13){$F_5$}
        \put(89,7){\scriptsize $-2$}
    \put(54,40){\line(1,0){30}}
    \put(92,40){\line(1,0){30}}
    \put(88,44){\line(0,1){30}}
    \put(88,36){\line(0,-1){30}}
\end{picture}
\end{center}
$w \geq 3. \quad$ If $w \geq 4$, then:
\[
  Z = F_1 + (F_2 + \cdots + F_5); \;\;
  \Delta = F_1 + \tfrac{1}{2}(F_2 + \cdots + F_5); \;\; \delta_y = 2.
\]
However, if $w = 3$, then $Z = 2F_1 + (F_2 + \cdots + F_5)$ (while $\Delta$
is the same), and $\delta_y = 1$. This exceptional case illustrates an
interesting property of the fundamental cycle; see (2.10) below.

\vspace{6pt}

{\bf (1.10)}
{\bf Exercise.} Calculate $Z$, $\Delta$, and $\delta_y$ for the following dual
graph (cf. \cite[p.350]{brieskorn}):

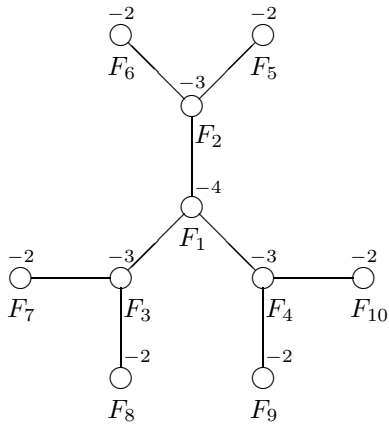
\begin{center}
\begin{picture}(180,180)(0,0)   
    \put(91,95){\scriptsize $-4$}
    \put(90,90){\circle{8}}
    \put(85,75){$F_1$}
       \put(90,94){\line(0,1){30}}
          \put(85,135){\scriptsize $-3$}
          \put(90,128){\circle{8}}
          \put(91,114){$F_2$}
              \put(92.82,130.82){\line(1,1){21.21}}
          \put(111.87,160.87){\scriptsize $-2$}
          \put(116.87,154.87){\circle{8}}
          \put(111.87,139.87){$F_5$}
              \put(87.18,130.82){\line(-1,1){21.21}}
          \put(58.13,160.87){\scriptsize $-2$}
          \put(63.13,154.87){\circle{8}}
          \put(58.13,139.87){$F_6$}
       \put(87.18,87.18){\line(-1,-1){21.21}}
          \put(58.13,69.13){\scriptsize $-3$}
          \put(63.13,63.13){\circle{8}}
          \put(64.13,48.13){$F_3$}
              \put(59.13,63.13){\line(-1,0){30}}
          \put(20.13,69.13){\scriptsize $-2$}
          \put(25.13,63.13){\circle{8}}
          \put(20.13,48.13){$F_7$}
              \put(63.13,59.13){\line(0,-1){30}}
          \put(64.13,31.13){\scriptsize $-2$}
          \put(63.13,25.13){\circle{8}}
          \put(58.13,10.13){$F_8$}
       \put(92.82,87.18){\line(1,-1){21.21}}
          \put(111.87,69.13){\scriptsize $-3$}
          \put(116.87,63.13){\circle{8}}
          \put(117.87,48.13){$F_4$}
              \put(116.87,59.13){\line(0,-1){30}}
          \put(117.87,31.13){\scriptsize $-2$}
          \put(116.87,25.13){\circle{8}}
          \put(111.87,10.13){$F_9$}
              \put(120.87,63.13){\line(1,0){30}}
          \put(149.87,69.13){\scriptsize $-2$}
          \put(154.87,63.13){\circle{8}}
          \put(149.87,48.13){$F_{10}$}
\end{picture}
\end{center}

\vspace{6pt}

{\bf (1.11)}
We conclude this section with another example. Assume that $f^{-1}(y)$
is a smooth curve $C$ of genus $g$, with $C^2=-w, w \geq 1$. This situation
can be realized easily in practice: for example, $y$ could be the vertex of
the cone $Y = \text{Proj} ( \bigoplus_{k \geq 0} H^0(C,kL) )$, with $C$ an
arbitrary smooth curve of genus $g$ and $L$ an arbitrary divisor of degree $w$
on $C$.

Then $Z=C, \;\; \Delta = ( \frac{2}{w}(g-1) + 1 ) C, \;\;
\text{ and } \;\; \delta_y = \dfrac{4(g-1)^2}{w}$.

If $g = 0$, then $\Delta = (1-\frac{2}{w})C$, $y$ is log-terminal, and
$\delta_y = \frac{4}{w}$. If $w=1$ then $y$ is smooth. If $w=2$ then $y$ is
an $A_1$ singularity (ordinary double point). If $w \geq 3$ then $y$ is a
log-terminal singularity ``of type $A_1$''.

If $g = 1$, then $y$ is log-canonical and $\delta_y = 0$. Such a singularity
is known as \emph{simply elliptic}.

If $g \geq 2$, then $y$ is not log-canonical. Note that $\delta_y$ can be
arbitrarily large in this case (if $g$ is large relative to $w$). In
particular, $\delta_y$ \emph{may} be greater than 4 (which is the value for
smooth points).


\section{A theorem of Reider type on normal surfaces with boundary}

\vspace{6pt}

In this section the ground field is \C.

\vspace{6pt}

{\bf (2.1)}
Let $Y$ be a projective surface over \C, and let $y$ be a fixed point on $y$.
Assume that $Y$ is smooth except perhaps at $y$, which may be either smooth or
a $RDP$. Let $M$ be a nef and big \Q-divisor on $Y$, with the property that 
$\rup{M}$ is Cartier.

Ein and Lazarsfeld proved the following criterion on global generation:

\vspace{4pt}

Let $B = \rup{M} - M$, and let $\mu = \mult_y B$ if $y$ is smooth, resp.
$\mu = \max \{ t \geq 0 \mid f^*B \geq tZ \}$ when $y$ is a $RDP$, where
$f:X\to(Y,y)$ is the minimal resolution and $Z$ is the fundamental cycle.

{\bf Theorem. (\cite[Theorem 2.3]{el})}
{\em
  Assume that $M^2 > (2-\mu)^2$ and $M \cdot \C \geq (2-\mu)$ for all curves
  $C$ through $y$ (when $y$ is smooth), resp. that $M^2 > 2 \cdot (1-\mu)^2$
  and $M \cdot C \geq (1-\mu)$ for all $C$ through $y$ (when $y$ is a $RDP$).

  Then $y \notin \bskm$.
}

\vspace{4pt}

This theorem was an important part of Ein and Lazarsfeld's proof of Fujita's
conjecture on smooth threefolds. In extending that work to (log-) teminal
threefolds (as required by the minimal model theory), it was necessary to
extend the criterion mentioned above to arbitrary normal surfaces. Such
extensions were obtained, e.g., in \cite[Theorem 1.6]{elm}, 
\cite[Theorem 7]{matsushita}. However, these generalizations, while
effective, are not optimal.

\vspace{4pt}

In a somewhat different direction, Kawachi and the author proved the following
criterion, of independent interest:

\vspace{4pt}

{\bf Theorem. (\cite[Theorem 2]{kawachim})}
{\em
  Let $Y$ be a normal surface, and let $y$ be a fixed point on $Y$. Let
  $\delta = \delta_y$ if $y$ is log-terminal, $\delta = 0$ otherwise. Let $M$
  be a nef divisor (with integer coefficients) on $Y$, such that $M^2>\delta$
  and $K_Y+M$ is Cartier.

  If $y \in \bskm$, then there exists an effective divisor $C$ passing through
  $y$, such that $M \cdot C < \frac{1}{2}\delta$ and $C^2 \geq M \cdot C -
    \frac{1}{4} \delta$
}
  (in particular, $y$ must be log-terminal, because $M$ is nef).

\vspace{4pt}

\begin{Remark}
  When $y$ is smooth, this is equivalent to Reider's original criterion
  (\cite[Theorem 1]{reider}). When $y$ is a $RDP$, we recover the
  Ein--Lazarsfeld criterion, plus a lower bound on $C^2$.
\end{Remark}

While this result has several applications to linear systems on normal
surfaces (cf. \cite{kawachim}), it cannot be used in the context of Fujita's
conjecture on terminal threefolds, because $M$ is required to  have integer
coefficients.

Kawachi formulated the following criterion for \Q-divisors on normal surfaces:

\vspace{4pt}
Let $Y,y$ be as before. Let $f:X \to (Y,y)$ be the minimal resolution if $y$
is singular, resp. the blowing-up at $y$ if $y$ is smooth. Let $Z$ and
$\Delta$ be the fundamental and the canonical cycle, respectively. Let
$\mu = \max \{ t \geq 0 \mid f^*B \geq t(Z-\Delta) \}$; note that
$\mu = 2 \cdot \mult_y Y$ if $y$ is smooth.

{\bf Open Problem. (cf. \cite{kawachi2}) }
{\em
  Let $\delta = \delta_y$ if $y$ is log-terminal, $\delta = 0$ otherwise.
  Let $M$ be a nef \Q-Weil divisor on $Y$, such that $K_Y+B+M$ is Cartier,
  where $B = \rup{M} - M$.

  If $M^2 > (1-\mu)^2 \delta$ and $M \cdot C \geq (1-\mu) \frac{\delta}{2}$
  for every curve $C$ through $y$, then $y \notin \bskm$.
}

\vspace{4pt}

When $M$ has integer coefficients, this criterion is the same as the one
mentioned above, minus the lower bound on $C^2$. Also, this criterion contains
the Ein--Lazarsfeld results for smooth and rational double points.

Kawachi formulated this as a theorem. Unfortunately his proof, based on a
case-by-case analysis, is incomplete. In this section we prove a slightly
weaker version, requiring that $M \cdot C \geq (1-\mu)$, rather than
$\geq (1-\mu) \frac{\delta}{2}$, for all $C$ through $y$. For application to
Fujita's conjecture on singular threefolds this makes little difference,
though, because $\delta_y$ cannot be controlled in that situation anyway; the
bound $\delta_y \leq 2$ (for $y$ singular) is used instead.

\vspace{6pt}

{\bf (2.2)}
Let $Y$ be a normal surface (= compact, normal two-dimensional algebraic space
over $\C$). Let $y \in Y$ be a given point, and let $B=\sum b_iC_i$ be an
effective $\Q$-Weil divisor on $Y$ with all $b_i \in \Q$, $0 \leq b_i \leq 1$;
the $C_i$ are distinct irreducible curves on $Y$. Since later we may need to
consider more curves $C_i$ than there are in $\Supp(B)$, we allow some 
coefficients $b_i$ to be $0$.  Let $f:X\to(Y,y)$ be the minimal resolution of
singularities of the germ $(Y,y)$ -- resp.  the blowing-up at $y$ if $y$ is a
smooth point. Let $f^{-1}(y) = \cup F_j, Z=\sum z_jF_j, \Delta=\sum a_jF_j,$
as in \S 1.

\vspace{6pt}

{\bf (2.3)}
Let $D_i = f^{-1}C_i$. Write $f^*B = f^{-1}B + B_{\text{exc}} =
	\sum b_i D_i + \sum b'_j F_j$.

\begin{definition}
  $(Y,B,y)$ is log-terminal (respectively log-canonical) if $a_j+b_j' < 1$
  (respectively $\leq 1$) for all $j$.

  Thus $(Y,y)$ is log-terminal (log-canonical) if $(Y,0,y)$ is.
\end{definition}

\begin{Remark}
  We do not require $K_Y+B$ to be $\Q$-Cartier at $y$ (unlike the similar
  definition in higher dimension); as in \S 1, this is a \emph{consequence} 
  of the other conditions. Note that $B \geq 0 \implies f^*B \geq 0$, and
  therefore $(Y,B,y)$ log-terminal (log-canonical) $\implies (Y,y)$
  log-terminal (log-canonical). Moreover, if $y \in \Supp(B)$, then
  $b'_j > 0,\; \forall j$, and therefore $(Y,B,y)$ log-canonical already
  implies $(Y,y)$ log-terminal.
\end{Remark}

\vspace{6pt}

{\bf (2.4)}
{\bf Definition.}
  Assume $(Y,B,y)$ is log-terminal. Define
  \[
      \mu=\mu_{B,y}=\max\{ t\geq 0 \mid f^*B \geq t(Z-\Delta) \}.
  \]

\begin{Remark}
  All the $z_j$ are $\geq 1$ and all the $a_j$ are
  $<1$; $\mu$ is given explicitly by
    \[
        \mu = \min \left\{ \frac{b_j'}{z_j-a_j} \right\} .
    \]
  Of course, $\mu = 0$ if $B=0$ (or, more generally, if $y \notin \Supp(B)$).
  Note that $\mu = \tfrac{1}{2} \mult_y(B)$ if $y$ is a smooth point of $Y$.
\end{Remark}

\begin{Lemma}
  If $(Y,B,y)$ is log-terminal and $\mu$ is defined as above, then
  $0\leq \mu <1$.
\end{Lemma}

Indeed, $\mu \geq 0$ is clear. On the other hand, if $\mu \geq 1$ then we have
$f^*B \geq (Z-\Delta)$, or $\Delta+f^*B \geq Z$, and therefore
$a_j + b_j' \geq z_j \geq 1$ for every $j$, contradicting log-terminality. \qed

\vspace{6pt}

{\bf (2.5)}
Let $Y$ be a normal surface and $y\in Y$ a given point, as in (2.1).  Let $M$
be a {\bf nef and big} $Q$-Weil divisor on $Y$ such that $K_Y+\rup{M}$ is
Cartier. Let $B = \rup{M} - M = \sum b_iC_i$; $B$ is an effective $Q$-Weil
divisor on $Y$, with $0 \leq b_i <1,\; \forall i$.

\vspace{4pt}

We will prove the following criteria for freeness at $y$:

\vspace{4pt}

\begin{Theorem}  \label{thm:notlt}
  If $(Y,B,y)$ is \textbf{not} log-terminal, then $y \notin \bskm$.
\end{Theorem}

\vspace{6pt}

In the following two theorems, assume that $(Y,B,y)$ \emph{is} log-terminal.
Then $(Y,y)$ is also log-terminal. Define $\mu$ as in (2.3) and $\delta_y$
as in \S 1.

Also, assume that $y$ is \emph{singular}.

\vspace{6pt}

\begin{Theorem}  \label{thm:kaw}
  Assume that $M^2>(1-\mu)^2\delta_y$ and $M \cdot C \geq (1-\mu)$ for every
  curve $C \subset Y$ passing through $y$. (Note that $M$ is still assumed to 
  be nef, i.e. $M \cdot C \geq 0$ even for curves $C \subset Y$ not passing 
  through $y$.)    						
\newline
  Then $y \notin \bskm$.
\end{Theorem}

\vspace{6pt}

In fact, if $y$ is a singularity of type $D_n$ or $E_n$ (see
\cite[Remark 9.7]{crepant}), we don't even need the assumption on $M \cdot C$
for $C$ through $y$:

\vspace{6pt}

\begin{Theorem}  \label{thm:dnen}
  Assume that $(Y,y)$ is a log-terminal singularity of type $D_n$ or $E_n$,
  and $M$ is nef and $M^2>(1-\mu)^2\delta_y$. If $(Y,y)$ is of type $D_n$,
  assume moreover that $M \cdot C > 0$ for every $C$ through $y$.
\newline
  Then $y \notin \bskm$.
\end{Theorem}

\vspace{6pt}

{\bf (2.6)}
First we reduce the proof of Theorems \ref{thm:notlt}, \ref{thm:kaw} and
\ref{thm:dnen} to the case when $y$ is the only singularity of $Y$:

\begin{Lemma}  \label{lemma:yonly}
  We may assume that $Y - \{y\}$ is smooth.
\end{Lemma}

\begin{proof}
  If $y$ is not the only singularity of $Y$, then let $g:S \to Y$ be a
  simultaneous resolution of all singularities of $Y$ \emph{except} $y$.
  Put $M' = g^*M$ and $y'=g^{-1}(y)$ (note that $g$ is an isomorphism of
  an open neighborhood of $y'$ onto an open neighborhood of $y$).
  $K_S+\rup{M'}$ is Cartier: outside $y'$ this is clear, because $S-\{y'\}$
  is smooth and $K_S+\rup{M'}$ has integer coefficients, and in a certain
  open neighborhood of $y'$ this is also clear, because $g$ is an isomorphism
  there, and $K_Y+\rup{M}$ is Cartier by hypothesis. Also, all the numerical
  conditions on $M$ are satisfied by $M'$ (in each of the hypotheses
  of Theorems \ref{thm:notlt}, \ref{thm:kaw} and \ref{thm:dnen}).
  If the theorems are true for $S$, $y'$, $M'$, then we get
  $y' \notin \Bs|K_S+\rup{M'}|$.

  Write $\Delta' = g^*K_Y-K_S$; $\Delta'$ is an effective $\Q$-divisor on $S$,
  and we have: $K_S+\rup{M'} = \ulc g^*K_Y-\Delta'+g^*\rup{M}-g^*B \urc =
  g^*(K_Y+ \rup{M}) - \rdn{\Delta' + g^*B}$ (note that $K_Y+\rup{M}$ is Cartier
  by hypothesis, and therefore $g^*(K_Y+\rup{M})$ has integer coefficients). 
  Write $N = \rdn{\Delta'+g^*B}$; $N$ is a divisor with integer coefficients
  on $S$, and $y' \notin \Supp(N)$, because in a certain neighborhood of $y'$,
  $\Delta'$ is zero and $g^*B$ is identified to $B$ -- whose coefficients
  are all $<1$.

  In the first part of the proof we found a section
  $s \in H^0(S, g^*(K_Y+\rup{M})-N)$ which doesn't vanish at $y'$.
  Multiplying $s$ by a global section of $\Ocal_S(N)$ whose zero locus is $N$,
  we find a new section $t \in H^0(S, g^*(K_Y+\rup{M}))$, which still 
  doesn't vanish at $y'$. In turn, $t$ corresponds to a global section of
  $\Ocal_Y(K_Y+\rup{M})$ which doesn't vanish at $y$.
\end{proof}

\vspace{6pt}

{\bf (2.7)}
Now assume that $y$ is the only singularity of $Y$. Let $f:X \to Y$ be the
minimal \emph{global} desingularization of $Y$,
$\Delta = f^*K_Y-K_X = \sum a_jF_j$, $B=\sum b_iC_i$, 
$f^*B=\sum b_iD_i + \sum b_j'F_j$, $Z=\sum z_j F_j$, etc.

\vspace{6pt}

\noindent {\bf Proof of Theorem \ref{thm:notlt}}

\vspace{3pt}

Assume that $(Y,B,y)$ is not log-terminal; then there is at least one $j$
such that $a_j+b_j' \geq 1$. $f^*M$ is nef and big on $X$, and therefore
the Kawamata--Viehweg vanishing theorem (cf. \cite[Lemma 1.1]{el}) gives:
$H^1(X, K_X+\rup{f^*M}) = 0$. 

$K_X+\rup{f^*M} = \ulc f^*K_Y-\Delta+f^*\rup{M}-f^*B \urc
  = f^*(K_Y+\rup{M})- \rdn{\Delta + f^*B}
  = f^*(K_Y+\rup{M})- \sum \rdn{a_j+b_j'} F_j
  = f^*(K_Y+\rup{M})- G$,
where $G$ is a nonzero effective divisor with integer coefficients on $X$
such that $f(G)=\{y\}$. ($G>0$ because $a_j+b_j'\geq 1$ for at least one $j$).

We have $H^1(X, f^*(K_Y+\rup{M})-G)=0$, and therefore the restriction map
\[
  H^0(X, f^*(K_Y+\rup{M})) \to H^0(G, f^*(K_Y+\rup{M}) |_G)
\]
is surjective.
As $f(G)=\{y\}$, $f^*(K_Y+\rup{M}) |_G$ is trivial, i.e. it has a global
section which doesn't vanish anywhere on $G$. By surjectivity, this section
lifts to a global section $s \in H^0(X,f^*(K_Y+\rup{M}))$ which doesn't
vanish anywhere on $G$. In turn, $s$ corresponds to a global section of
$\Ocal_Y(K_Y+\rup{M})$ which doesn't vanish at $y$, or else $s$ would
vanish everywhere on $f^{-1}(y)$, and in particular on $G$.   	\qed

\vspace{6pt}

{\bf (2.8)}
Assume that $y$ is log-terminal but not smooth.
Thus $\Delta=\sum a_j F_j$ with $0\leq a_j <1$ for every $j$.

\vspace{6pt}
\noindent {\bf Proof of Theorem \ref{thm:kaw}}

\vspace{3pt}

\begin{Lemma}
  If $M^2 >(1-\mu)^2\delta_y$, then we can find an effective $\Q$-Weil divisor
  $D$ on $Y$ such that $D\equiv M$ and $f^*D\geq(1-\mu)(Z-\Delta)$.
\end{Lemma}

\begin{proof}
  Since $M^2>(1-\mu)^2\delta_y$, $f^*M-(1-\mu)(Z-\Delta)$ is in the positive
  cone of $X$ (see \cite{kawachim}, (2.3) for a similar argument). In
  particular, $f^*M-(1-\mu)(Z-\Delta)$ is big.

  Let $G\in \big|\,k(f^*M-(1-\mu)(Z-\Delta)\,)\,\big|$ for some $k$
  sufficiently large and divisible. Let $T=\frac{1}{k}\,G+(1-\mu)(Z-\Delta)$.
  Write $T=\sum d_i D_i+ \sum t_j F_j$. Define $D=f_*T=\sum d_iC_i$, and write
  $f^*D=\sum d_i D_i+ \sum d_j'F_j$. We have:
  \begin{enumerate}
    \item $T \equiv f^*M$, because $kT \sim kf^*M$;
    \item $T\geq 0$, and therefore $D \geq 0$;
    \item $T\cdot F_j=0$ for every exceptional curve $F_j$, because
			$T\equiv f^*M$;
    \item $f^*D=T$; indeed, the coefficients $d_j'$ are uniquely
  			determined by the condition $f^*D \cdot F_j = 0$ for
			every $j$, and $T$ already satisfies this condition;
    \item $f^*D=T\geq(1-\mu)(Z-\Delta)$, because $G\geq 0$;
    \item finally, $D\equiv M$, because $f^*D=T\equiv f^*M$.
  \end{enumerate}
\end{proof}
\begin{Remark}
  $f^*D\geq(1-\mu)(Z-\Delta)$ means $d_j'\geq(1-\mu)(z_j-a_j)$
  for every $j$. We may assume, however, that $d_j' > (1-\mu)(z_j-a_j)$ for
  every $j$. Indeed, since $M^2>(1-\mu)^2\delta_y$, we have 
  $M^2>(1-\mu)^2\delta_y(1+\epsilon)^2$ for some small rational number
  $\epsilon>0$; then working as above we can find $D \equiv M$ such that
  $d_j'\geq(1-\mu)(z_j-a_j)(1+\epsilon)>(1-\mu)(z_j-a_j)$ for every $j$.
\end{Remark}

We'll assume that {\boldmath $d_j'>(1-\mu)(z_j-a_j)$} for all $j$.

\vspace{8pt}

For every rational number $c$, $0<c<1$, let $R_c=f^*(M-cD)$.
$R_c\equiv(1-c)f^*M$, so that $R_c$ is nef and big, and we have
\begin{equation}
  H^1(X,K_X+\rup{R_c})=0.		\label{eq:van}
\end{equation}
\begin{align*}
    K_X+\rup{R_c}     &= \rup{f^*K_Y-\Delta+f^*\rup{M}-f^*B-cf^*D}  \\
		      &= f^*(K_Y+\rup{M})-\rdn{\Delta+f^*B+cf^*D}    \\
      &= f^*(K_Y+\rup{M})-\sum\rdn{b_i+cd_i}D_i-\sum\rdn{a_j+b_j'+cd_j'}F_j.
\end{align*}

Choose $c=\min\left\{\frac{1-a_j-b_j'}{d_j'}, \text{ all $j$};
  \frac{1-b_i}{d_i}, \text{ all $i$ such that $d_i>0$ and $y\in C_i$}\right\}$.

Note that $c > 0$, because $(Y,B,y)$ is log-terminal, and $c<1$, because
for every $j$ we have $b_j'\geq\mu(z_j-a_j)$ and $z_j\geq 1$, and therefore
$1-a_j-b_j' \leq z_j-a_j-\mu(z_j-a_j) = (1-\mu)(z_j-a_j) < d_j'$. Therefore
(\ref{eq:van}) holds for this choice of $c$.

Note also that $0<a_j+b_j'+cd_j'\leq 1$ for all $j$; so
$F \eqdef \sum \rdn{a_j+b_j'+cd_j'} F_j$ is either zero or a sum of distinct
irreducible components, $F=F_1+\cdots+F_s$ (after re-indexing the $F_j$ if
necessary). Similarly, $\sum \rdn{b_i+cd_i}D_i = N+A$, where
$\Supp(N) \cap f^{-1}(y)=\emptyset$, and $A$ is either zero or a sum of
distinct irreducible components, $A=D_1+\cdots+D_t$, where $D_1, \ldots, D_t$
meet $f^{-1}(y)$. Each component $F_j$ of $F$ (if any), and each component $D_i$
of $A$ (if any), has coefficient $1$ in $\Delta + f^*B + c f^*D$.  Also, $F$
and $A$ cannot both be equal to zero.

\vspace{8pt}

We will use the following form of the Kawamata--Viehweg vanishing theorem
(see, for example, \cite{el}, Lemma 2.4):

\begin{Lemma} 		 \label{lemma:vt2} 
  Let $X$ be a smooth projective surface over $\C$, and let $R$ be a nef
  and big $\Q$-divisor on $X$. Let $E_1,\ldots,E_m$ be distinct irreducible
  curves such that $\rup{R} \cdot E_i > 0$ for every $i$. Then
  $$
	H^1(X,K_X+\rup{R}+E_1+\cdots+E_m)=0.
  $$
\end{Lemma}

\vspace{8pt}

We consider two cases, according to whether $F \neq 0$ or $F=0$.

\vspace{6pt}

\textbf{Case I:} $\mathbf{F \neq 0.}$

Using Lemma \ref{lemma:vt2} for $R_c$ in place of $R$ and $D_1,\ldots,D_t$
in place of $E_1,\ldots,E_m$, we get $H^1(X,K_X+\rup{R_c}+A)=0$, or
$H^1(X,f^*(K_Y+\rup{M})-N-F)=0$. We conclude as in the proof of
Theorem \ref{thm:notlt}; the `` $-N$ '' part is treated as in the proof
of Lemma \ref{lemma:yonly}.
\emph{Note:} $R_c \cdot D_i = (1-c) M \cdot C_i > 0$ for every $i$,
and $D_i$ has integer coefficient in $R_c$ if $D_i$ is a component of $A$ --
and therefore $\rup{R_c} \cdot D_i > 0$ for such $D_i$.

\vspace{6pt}

\textbf{Case II:} $\mathbf{F = 0.}$

As noted earlier, in this case $A \neq 0$; using Lemma \ref{lemma:vt2} as in
Case I above, we get $H^1(X,f^*(K_Y+\rup{M})-N-D_1)=0$, and therefore 
the restriction map
\begin{equation}
  H^0(X,f^*(K_Y+\rup{M})-N) \to H^0(D_1,(f^*(K_Y+\rup{M})-N|_{D_1})
							\label{eq:surj}
\end{equation}
is surjective.

$D_1 \cap f^{-1}(y) \neq \emptyset$; let $x\in D_1\cap f^{-1}(y)$.
Assume we can find a section $s'\in H^0(D_1, f^*(K_Y+\rup{M})-N|_{D_1})$
such that $s'(x)\neq 0$. Then by the surjectivity of (\ref{eq:surj}) we can 
find $s\in H^0(X,f^*(K_Y+\rup{M})-N)$ such that $s(x)\neq 0$; then we
conclude as in the proof of Lemma \ref{lemma:yonly}.

Hence the proof is complete if we show that $x \notin \Bs|\,
f^*(K_Y+\rup{M})-N|_{D_1}\,|$. Note that $f^*(K_Y+\rup{M})-N|_{D_1}=
K_X+\rup{R_c}+A|_{D_1}=(K_X+D_1)|_{D_1}+(\rup{R_c}+D_2+\cdots+D_t)
|_{D_1}=K_{D_1}+(\rup{R_c}+D_2+\cdots+D_t)|_{D_1}$. By \cite{har},
Theorem 1.4 and Proposition 1.5, it suffices to show that $(\rup{R_c}+
D_2+\cdots+D_t)\cdot D_1 >1$ (then this intersection number is $\geq 2$,
because it is an integer).

$\rup{R_c} = R_c+\sum\{b_i+cd_i\}D_i+\sum\{a_j+b_j'+cd_j'\}F_j$. Note that
$b_1 + c d_1 = 1$, and therefore $\{b_1+cd_1\}=0$. Also, since $F=0$, we have
$0\leq a_j+b_j'+cd_j'<1$ for every $j$ -- and therefore
$\{a_j+b_j'+cd_j'\}=a_j+b_j'+cd_j'$. Hence we have
$$
  (\rup{R_c}+D_2+\cdots+D_t)\cdot D_1 \geq R_c\cdot D_1 +
		\sum (a_j+b_j'+cd_j')F_j \cdot D_1.
$$
$R_c \equiv (1-c)f^*M$, so that $R_c\cdot D_1 = (1-c)M\cdot C_1 \geq
(1-c)(1-\mu)$, because $y \in C_1 = f_*D_1$. Also, $D_1$ meets at least
one $F_j$, say $F_1$. Therefore we have
\begin{align*}
  (\rup{R_c}+D_2+\cdots+D_t)\cdot D_1 &\geq(1-c)(1-\mu)+(a_1+b_1'+cd_1') \\
  	&> (1-c)(1-\mu) +a_1 +\mu(z_1-a_1) +c(1-\mu)(z_1-a_1)                \\
	&\geq (1-c)(1-\mu) +a_1 +\mu(1-a_1) +c(1-\mu)(1-a_1)                 \\
	&= 1 +(1-c)(1-\mu)a_1						     \\
	&\geq 1.
\end{align*} 	\qed 			

\vspace{6pt}

{\bf (2.9)}
Finally, we prove Theorem \ref{thm:dnen}. We assume the reader is familiar
with the classification of $RDP$'s. The classification of log-terminal
singularities is similar: if $f:X\to (Y,y)$ is the minimal resolution of a
log-terminal germ and $f^{-1}(y)=\cup F_j$, then the $F_j$ are smooth rational
curves, and the dual graph is a graph of type $A_n$, $D_n$ or $E_n$ (see,
e.g., \cite[\S 9]{crepant}, or \cite{alexeev}).  The only difference is that
in the log-terminal case the self-intersection numbers $-w_j = F_j^2$ are not
necessarily all equal to $-2$. 

The classification of log-terminal singularities with non-zero reduced boundary
is even simpler. We have:

\begin{Lemma}   \label{lemma:alexeev}
  Assume $(Y,y)$ is a normal surface germ, and $C_1$ is a reduced, irreducible
  curve on $Y$ such that $y \in C_1$. If $(Y,C_1,y)$ is log-terminal, then
  $(Y,y)$ is of type $A_n$. More explicitly, if $f:X\to (Y,y)$ is the minimal
  resolution, $f^{-1}(y) = F_1 \cup \ldots \cup F_n$, $-w_j=F_j^2$, and
  $f^{-1}C_1=D_1$, then the dual graph of the resolution is:
    
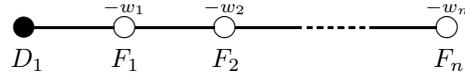
\begin{center}
    \begin{picture}(190,60)(0,0)  
        \put(14,30){\circle*{8}}
        \put(9,15){$D_1$}
            \put(18,30){\line(1,0){30}}
        \put(44,36){\scriptsize $-w_1$}
        \put(52,30){\circle{8}}
        \put(47,15){$F_1$}
            \put(56,30){\line(1,0){30}}
        \put(82,36){\scriptsize $-w_2$}
        \put(90,30){\circle{8}}
        \put(85,15){$F_2$}
            \put(94,30){\line(1,0){25}}
            \multiput(121,30)(4,0){6}{\line(1,0){2}}
            \put(145,30){\line(1,0){25}}
        \put(166,36){\scriptsize $-w_n$}
        \put(174,30){\circle{8}}
        \put(169,15){$F_n$}
    \end{picture}
    \end{center}
  (If $y$ is smooth, then $n=1$ and $w_1=1$; otherwise all $w_j\geq 2$.)
\end{Lemma}

Similarly, for log-canonical singularities with boundary we have:

\begin{Lemma} 	\label{lemma:typeen}
  Let $(Y,y)$ and $C_1$ be as in the previous lemma. If $(Y,C_1,y)$ is
  log-canonical, then either $(Y,y)$ is of type $A_n$ and the dual graph
  is the one shown above, or $(Y,y)$ is of type $D_n$ and the dual graph
  of the resolution is:
    
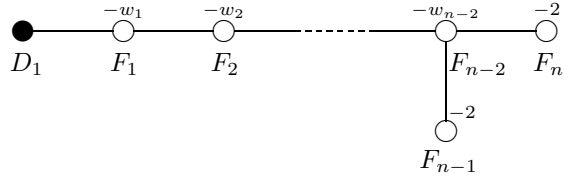
\begin{center}
    \begin{picture}(230,98)(0,-38)   
        \put(14,30){\circle*{8}}
        \put(9,15){$D_1$}
            \put(18,30){\line(1,0){30}}
        \put(44,36){\scriptsize $-w_1$}
        \put(52,30){\circle{8}}
        \put(47,15){$F_1$}
            \put(56,30){\line(1,0){30}}
        \put(82,36){\scriptsize $-w_2$}
        \put(90,30){\circle{8}}
        \put(85,15){$F_2$}
            \put(94,30){\line(1,0){25}}
            \multiput(121,30)(4,0){6}{\line(1,0){2}}
            \put(145,30){\line(1,0){25}}
        \put(161,36){\scriptsize $-w_{n-2}$}
        \put(174,30){\circle{8}}
        \put(175,15){$F_{n-2}$}
        \put(178,30){\line(1,0){30}}
        \put(207,36){\scriptsize $-2$}
        \put(212,30){\circle{8}}
        \put(207,15){$F_n$}
        \put(174,26){\line(0,-1){30}}
        \put(175,-3){\scriptsize $-2$}
        \put(174,-8){\circle{8}}
        \put(164,-22){$F_{n-1}$}
    \end{picture}
    \end{center}
\end{Lemma}

All these facts can be found, for example, in \cite{alexeev}, and also in
\cite{crepant}.

\vspace{8pt}

\noindent {\bf Proof of Theorem \ref{thm:dnen}}

\vspace{3pt}

Going back to the proof of theorem \ref{thm:kaw}, we will show that Case II of
the proof is not possible if $y$ is of type $D_n$ or $E_n$. Therefore the
condition $M \cdot C \geq (1-\mu)$ for $C$ through $y$ is no longer necessary.
We still need $M \cdot C > 0$ for $C$ through $y$ if $A \neq 0$, to be able to
use Lemma \ref{lemma:vt2} (see the \emph{Note} at the end of Case I). However,
we will show that $A = 0$ if $y$ is of type $E_n$, so that in that case we
don't even need the condition $M \cdot C > 0$ for $C$ through $y$.

Assume that $A \neq 0$. Let $C_1$ be a component of $A$, and let
$f^*C_1 = D_1 + \sum c'_j F_j$. Since $C_1$ is a component of $A$, we have
$b_1+cd_1=1$, and therefore $B+cD=C_1+\textit{ other terms}\,$; consequently
$f^*B+cf^*D \geq f^*C_1$, and in particular $b'_j+cd'_j \geq c'_j, \forall j$.

If we end up in Case II in the proof of Theorem \ref{thm:kaw}, then
$a_j+b'_j+cd'_j < 1$ for all $j$; consequently $a_j+c'_j < 1, \forall j$, so
that $(Y,C_1,y)$ is log-terminal. By Lemma \ref{lemma:alexeev}, $(Y,y)$ must
be of type $A_n$.

Now assume that we end up in Case I, with $A \neq 0$. We still have
$a_j+b'_j+cd'_j \leq 1, \forall j$, and therefore $(Y,C_1,y)$ is log-canonical.
By Lemma \ref{lemma:typeen}, $(Y,y)$ is either of type $A_n$ or of type $D_n$.
		\qed

\begin{Remark}
For singularities of type $D_n$ and $E_n$, Theorem \ref{thm:dnen} is stronger
than the Open Problem (see (2.1)). To complete the proof of the Open Problem, 
the only case to consider is that of a singularity of type $A_n$.

In the proof of Theorem \ref{thm:kaw}, Case II (which is possible only when
$y$ is of type $A_n$), we used the inequality $M \cdot C \geq (1-\mu)$ for
$y \in C$. In fact, using Lemma \ref{lemma:alexeev} and modifying slightly the
final computation in (2.8), we see that we need slightly less: $M \cdot C \geq
(1-\mu)(1-a)$, where $a = \min \{ a_1,a_n \}$ (note that in Lemma
\ref{lemma:alexeev}, $D_1$ could meet either $F_1$ or $F_n$).

On the other hand, for $y$ of type $A_n$ we have $\delta_y = 2-(a_1+a_n)$
(for $n=1$ this follows from (1.11); for $n \geq 2$ use Lemma 2, (b), and the
obvious formulae $(Z-\Delta) \cdot F_j = -1$ for $j=1$ and $j=n$, 
$(Z-\Delta) \cdot F_j = 0$ otherwise). Thus the Open Problem requires that 
$M \cdot C \geq (1-\mu)(1- \frac{a_1+a_n}{2})$ for $y \in C$. In particular,
the Open Problem is proved if $a_1=a_n$ (e.g., if $y$ is of type $A_1$).
\end{Remark}

\vspace{6pt}

{\bf (2.10)}
Analyzing the proofs of Theorems \ref{thm:kaw} and \ref{thm:dnen}, we may ask:
what was the relevance of $Z$ being the fundamental cycle of $y$? $\Delta$
arises naturally, as $f^*K_Y-K_X$; but $Z$ could have been any effective, 
$f$-exceptional cycle with integer coefficients such that $z_j \geq 1$ for all
$j$ (i.e., such that $\Supp(Z) = f^{-1}(y)$). The answer is provided by the
following proposition:

\begin{Proposition}	\label{prop:caract}
  Let $(Y,y)$ be a log-terminal singularity, with $f:X \to (Y,y)$ the minimal
  resolution (resp. the blowing-up at $y$ if $y$ is smooth). Let $Z$ and
  $\Delta$ be the fundamental, resp. the canonical cycle. Let $Z'=\sum z'_jF_j$
  be any other effective, $f$-exceptional cycle (with integer coefficients),
  such that $z'_j \geq 1$ for all $j$.

  Then $\delta_y \leq \delta'$, where $\delta_y = -(Z-\Delta)^2$ and $\delta' =
  -(Z'-\Delta)^2$. Moreover, $Z$ is the (unique) largest cycle among all the
  $Z'$ for which $\delta' = \delta_y$.
\end{Proposition}
The proof depends on the detailed classification of log-terminal surface
singularities, cf. \cite{crepant}. Explicitly, we need the following lemma,
which can be proved by brute force (the computations are straightforward in
all cases; for reference, they can be found in \cite{kawachi1}):

\begin{Lemma}	\label{lemma:tech}
  Let $(Y,y)$ be log-terminal.

  (a) There is at most one $F_j$ with $(Z-\Delta) \cdot F_j > 0$. If such an
  $F_j$ exists, then $(Z-\Delta) \cdot F_j = 1$ and the corresponding
  $w_j = -F_j^2 \geq 3$.

  (b) There is at most one $F_j$ with $z_j \geq 2$ and $(Z-\Delta)\cdot F_j<0$.
  If such an $F_j$ exists, then $z_j=2$ and $(Z-\Delta) \cdot F_j = -1$.
\end{Lemma}

\vspace{8pt}

\noindent {\bf Proof of Proposition \ref{prop:caract}}

\vspace{3pt}

Let $Z' = Z+(P-N)$ with $P,N \geq 0$ without common components. Then
\begin{align*}
  \delta' &= -(Z'-\Delta)^2 = -(Z+P-N-\Delta)^2				\\
  	  &= -(Z-\Delta)^2 - (P-N)^2 - 2(Z-\Delta) \cdot (P-N)    	\\
  	  &= \delta_y + (-P^2) - 2 (Z-\Delta) \cdot P + (-N^2) +
  	  		2 (Z-\Delta) \cdot N + 2 (P \cdot N).
\end{align*}
Since $P \cdot N \geq 0$, the Proposition is proved if we can prove that
\begin{enumerate}
  \item[(a)] $(-P^2) - 2 (Z-\Delta) \cdot P > 0$ if $P>0$;
  \item[(b)] $(-N^2) + 2 (Z-\Delta) \cdot N \geq 0$.
\end{enumerate}

\vspace{5pt}

\noindent \emph{Proof of (a).}

\vspace{3pt}

If $P = \sum t_jF_j$, then we have
\[
    (-P^2) - 2 (Z-\Delta) \cdot P = (-P^2) - 2\sum t_j (Z-\Delta) \cdot F_j.
\]
If $(Z-\Delta) \cdot F_j \leq 0$ for all $j$, then we are done (note that
$(-P^2) > 0$ if $P>0$). Otherwise there is exactly one $j$, call it $j_0$,
such that $(Z-\Delta) \cdot F_{j_0} > 0$. By Lemma \ref{lemma:tech}, we have
$w_{j_0} = -F_{j_0}^2 \geq 3$ and $(Z-\Delta) \cdot F_{j_0} = 1$. Therefore
\[
   (-P^2) - 2\sum t_j (Z-\Delta) \cdot F_j \geq (-P^2) - 2t_{j_0}.
\]
We will show that $(-P^2) \geq t_{j_0}^2 + 2$; then (a) will follow.

Write $F_j^2 = -w_j, F_i \cdot F_j = l_{ij}$ for $i<j$ ($l_{ij} = 1$ if $F_i$
meets $F_j$, $0$ otherwise). We have:
\begin{align*}
  (-P^2) &= \sum_j w_j t_j^2 - \sum_{i<j} 2l_{ij} t_i t_j	\\
  	 &\geq t_{j_0}^2 + \sum_j 2t_j^2 - \sum_{i<j} 2 l_{ij} t_i t_j
\end{align*}
(note that $w_j \geq 2, \, \forall j$, and $w_{j_0} \geq 3$).

Now consider a singularity $(Y',y')$, whose minimal resolution has a dual
graph identical to that of $(Y,y)$, except that ${F'_j}^2 = -2,\,\forall j$
(thus $y'$ is a ``true'' $A_n, D_n,$ or $E_n$ rational double point).
If $P' = \sum t_j F'_j$ (having the same coefficients as $P$), then
$(-{P'}^2) > 0$; that is,
\[
    \sum_j 2t_j^2 - \sum_{i<j} 2l_{ij} t_i t_j > 0
\]
--- and therefore $\geq 2$, because it is an \emph{even} integer.   \qed

\vspace{5pt}

\noindent \emph{Proof of (b).}

\vspace{3pt}

If $N = \sum x_j F_j$, then we have
\[
    (-N^2) + 2 (Z-\Delta) \cdot N = (-N^2) + 2 \sum x_j (Z-\Delta) \cdot F_j.
\]
If $(Z-\Delta) \cdot F_j \geq 0$ for all $j$, or if $x_j = 0$ whenever
$(Z-\Delta) \cdot F_j < 0$, then we are done. Note that $x_j = z_j - z'_j$ if
$z_j > z'_j$, and $0$ otherwise. Thus $x_j \geq 1 \implies z_j \geq 2$.
Therefore, by Lemma \ref{lemma:tech}, there can be at most one negative term
in $2\sum x_j (Z-\Delta)\cdot F_j$; and if there is one, corresponding, say,
to $j_1$, then $z_{j_1} = 2$ and $(Z-\Delta) \cdot F_{j_1} = -1$.
Therefore $x_{j_1}=1$, and $2(Z-\Delta)\cdot N \geq -2$. Finally, 
$(-N^2) \geq 2$ (if $N \neq 0$), as in the proof of (a) above.    \qed

\vspace{5pt}

\noindent \emph{Remarks.}

\emph{1.} We showed that $\left[ -(Z'-\Delta)^2 = \delta_y \right] \implies 
[Z' \leq Z]$. The converse is not always true. For example, if $F_j$ is a
component with $z_j=2$ and $(Z-\Delta)\cdot F_j = 0$ (such components exist in
some cases of type $D_n$ and $E_n$, cf. \cite{kawachi1}), then taking
$Z'=Z-F_j$ we get $-(Z'-\Delta)^2 > \delta_y$ (cf. the proof of (b) above).

\emph{2.} If $(Y,y)$ is not log-terminal, then the statement of Proposition
\ref{prop:caract} is no longer necessarily true, even if $y$ is rational.
For example:
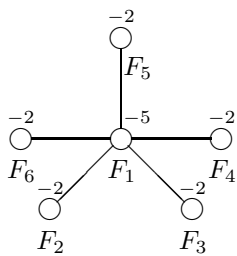
\begin{center}
\begin{picture}(100,110)(0,0)
    \put(51,56){\scriptsize $-5$}
    \put(50,50){\circle{8}}
    \put(45,35){$F_1$}
    	\put(47.17,47.17){\line(-1,-1){21.21}}
    \put(18.13,29.13){\scriptsize $-2$}
    \put(23.13,23.13){\circle{8}}
    \put(18.13,8.13){$F_2$}
    	\put(52.83,47.17){\line(1,-1){21.21}}
    \put(71.87,29.13){\scriptsize $-2$}
    \put(76.87,23.13){\circle{8}}
    \put(71.87,8.13){$F_3$}
    	\put(54,50){\line(1,0){30}}
    \put(83,56){\scriptsize $-2$}
    \put(88,50){\circle{8}}
    \put(83,35){$F_4$}
    	\put(50,54){\line(0,1){30}}
    \put(45,94){\scriptsize $-2$}
    \put(50,88){\circle{8}}
    \put(51,74){$F_5$}
    	\put(46,50){\line(-1,0){30}}
    \put(7,56){\scriptsize $-2$}
    \put(12,50){\circle{8}}
    \put(7,35){$F_6$}
\end{picture}
\end{center}
$Z = F_1 + (F_2 + \cdots + F_6), \; \Delta = \frac{6}{5} F_1 + 
\frac{3}{5} (F_2 + \cdots + F_6)$, and $\delta_y = \frac{13}{5}$;     \newline
but $\delta' = -(Z'-\Delta)^2 = \frac{8}{5} < \delta_y$ for $Z'=Z+F_1$.

\vspace{5pt}

\noindent \emph{Exercise.} Is the statement of Proposition \ref{prop:caract}
true for log-canonical singularities?

Hint: There is nothing to prove if $\delta_y \leq 1$. Notice that in the
exceptional case of (1.9), Type 3 with $w=3$, we have $\delta' = 2$ for
$Z' = F_1 + (F_2 + \cdots + F_5)$, while $\delta_y=1$ (in that case we have
$Z = 2 F_1 + (F_2 + \cdots + F_5)$).





\end{document}